\documentclass[11pt, twoside]{amsart}
\usepackage{subfigure}
\usepackage{epsfig}

\usepackage{latexsym}
\usepackage{graphicx}
\usepackage[centertags]{amsmath}
\usepackage{amsfonts}
\usepackage{amssymb}
\usepackage{amsthm}
\usepackage{newlfont}
\usepackage{subfigure}
\usepackage{epsfig}

\theoremstyle{plain}
\newtheorem{theorem}{Theorem}[section]

\newtheorem{corollary}[theorem]{Corollary}
\newtheorem{lemma}[theorem]{Lemma}

\theoremstyle{definition}
\newtheorem{definition}[theorem]{Definition}
\newtheorem{remark}[theorem]{Remark}

\newcommand{\E}[1][3]{\mathbb{E}^{#1}}
\newcommand{\R}[1][3]{\mathbb{R}^{#1}}
\newcommand{\s}[1][1]{\mathbf{S}^{#1}}
\newcommand{\h}[1][3]{\mathbb{H}^{#1}}

\newcommand{\cross}{\times}
\newcommand{\bndry}{\partial}

\newcommand{\pii}{\pi_1}

\newcommand{\tld}{\tilde}

\newcommand{\tc}{\cal T}
\newcommand{\adj}{adj_{\mathbf B_{r,v}^s}}

\newcommand{\homeo}{\cong}

\begin{document}

\title[Universal cover of 3-manifolds built from injective handlebodies]{The universal cover of 3-manifolds built from injective handlebodies is $\R[3]$.}

\author{James Coffey}
\thanks{Research partially supported by the Australian Research Council.}
\address{Department of Mathematics and Statistics. University of Melbourne, Australia, 3010}
 \email{coffey@ms.unimelb.edu.au}

\begin{abstract}
This paper gives a proof that the universal cover of a closed 3-manifold built from three $\pi_1$-injective handlebodies is homeomorphic to $\mathbb R^3$.  This construction is an extension to handlebodies of the conditions for gluing of three solid tori to produce non-Haken Seifert fibered manifolds with infinite fundamental group. This class of manifolds has been shown to contain non-Haken non-Seifert fibered manifolds.
\end{abstract}

\maketitle

\section{Introduction}

It is a long standing conjecture that all closed $P^2$-irreducible 3-manifolds with infinite fundamental group have universal cover homeomorphic to $\R[3]$.  In 1969 this was shown to be true in the case that the manifold is Haken by F.\,Waldhausen in \cite{wa1}. Herbert Seifert showed, in his 1932 doctoral dissertation ``Topology of 3-dimensional fibred spaces'',   that the conjecture is true for Seifert fibered manifolds with infinite fundamental group.   In 1987 J.\,Hass, H.\,Rubinstein and P.\,Scott, in \cite{hrs1}, showed that the conjecture was true if the manifold is $P^2$-irreducible and it contains an essential surface other than $S^2$ or $P^2$.  If Thurston's Geometrisation conjecture is proven, which looks likely, following the announcements of Perelman, then the conjecture will be proven as a result.  However it is still an interesting question to answer, using more direct means, for topological classes of manifolds that do not appear to have any obvious unifying underlying geometry. Note also that a similar question arises for higher dimensional aspherical manifolds and the methods in this paper may extend to suitable classes of examples.

The class of manifolds being considered in this paper is said to meet the `disk-condition'.  The construction of these manifolds has been discussed by H.\,Rubinstein and the author in \cite{C&R1}.  This construction is an extension to handlebodies of the gluing of three solid tori which produces non-Haken Seifert fibered manifolds with infinite fundamental group. A definition is given in the next section. The structure given by this construction tells us a lot about this class of  manifolds. In \cite{C&R1} it is shown that the characteristic variety of these manifolds can be constructed in a rather elegant way. Also examples of atoroidal manifolds and non-Haken non-Seifert fibered manifolds that meet the disk-condition are given. In \cite{Cof1} the author gave an algorithm to show that the word problem in the fundamental group of manifolds that meet the disk-condition is solvable. The main theorem proven in this paper is:

\begin{theorem}\label{theorem: universal cover is $R^3$}
All manifolds that meet the disk-condition have a universal cover homeomorphic to $R^3$.
\end{theorem}

In 1961 M.\,Brown in \cite{br} showed that if an $n$ dimensional metric space contained an infinite sequence of expanding balls, such that the union of the sequence is the space then the space is homeomorphic to $\R[n]$. So the general idea for the proof of this theorem is to produce an infinite expanding sequence of open balls that meets Brown's conditions, in the universal cover of manifolds that meet the disk-condition.  This is done by taking the dual 2-complex to the structure in the universal cover and then carefully defining an expanding sequence of simplicial subsets, such that the interior of the corresponding manifold in the universal cover is an open ball. This approach is quite pleasing  as the disk-condition is used implicitly  by the observation that the dual 2-complex turns out to be a CAT$(0)$ metric space.

\section{Preliminaries and definitions}
\label{section: universal cover-preliminaries}

Throughout this paper we will assume that, unless stated otherwise, we are working in the PL category of manifolds and maps.  Even though we will not explicitly use this structure we will use ideas that are a consequence, such as regular neighbourhoods and transversality as defined by C.\,Rourke and B.\,Sanderson in \cite{rou&sa1}. The standard definitions in this field, as given by J.\,Hempel in \cite{Hem1} or W. Jaco in \cite{Ja1}, are used.

A manifold $M$ is \textbf{closed} if $\bndry M = \emptyset$ and \textbf{irreducible} if every embedded $S^2$ bounds a ball.  We will assume that, unless otherwise stated, all 3-manifolds are orientable. The reason for this is that all closed non-orientable $\mathbb P^2$-irreducible 3-manifolds are Haken. (A manifold is $\mathbb P^2$-irreducible  if it is irreducible and does not contain any embedded $2$-sided projective planes).

If $M$ is a 3-manifold and $S$ a surface, which is not a sphere, disk or projective plane,  then a map $f:S\to M$ is called \textbf{essential} if the induced map $f^*:\pii(S)\to \pii(M)$ is injective. This is also known as a $\mathbf{\pi_1}$\textbf{-injective} map.   Also $f:S \to M$ is a \textbf{proper map} if $f(S)\cap \bndry M  = f(\bndry S)$.  Unless otherwise stated,  a homotopy/isotopy of a proper map, is assumed to be proper. That is, at each point the map remains proper.  To reduce  notation an isotopy/homotopy of a surface $S\subset M$ are used with out defining the map.  Here we are assuming that there is a map $f:S \to M$ and we are referring to  an isotopy/homotopy of $f$, however defining the map is often unnecessary and would only add to excessive book keeping.

If $H$ is a handlebody and $D$ is a properly embedded disk in $H$ such that $\bndry D$ is essential in $\bndry H$ then $D$ is a  \textbf{meridian disk} of $H$.  If $D$ is a proper singular disk in $H$ such that $\bndry D$ is essential in $\bndry H$, then it is called a \textbf{singular meridian disk}.

\begin{definition}
    For $H$ a handlebody, $\tc$ a set of curves in $\bndry H$ and $D$ a meridian disk, let $|D|$ be the number of intersection between $D$ and $\tc$.
\end{definition}

By the result of Freedman, Hass and Scott in \cite{fhs1}, if we put a metric on $\bndry H$ and isotop both $\tc$ and $D$ so that $\tc$ and $\bndry D$ are length minimising in $\bndry H$ then the number of intersections will be minimal. Note that when there are parallel curves in $\tc$ we need to `flatten' the metric in their neighbourhood so they remain disjoint.  Let $\mathcal D$  be a set containing a single representative from each isotopy class of meridian disks of $H$. We can also assume that the boundary of each meridian disk in $\mathcal D$ is a length minimising geodesic.  Therefore the number of intersections between the boundaries of two disks in $\cal D$ is minimal. In  \cite{C&R1} there is the following lemma proven.

\begin{lemma}\label{lemma: m-disk intersection.}
    Any two disks of $\mathcal D$ can be isotoped, leaving their boundaries fixed, so that the curves of intersections are  properly embedded arcs.
\end{lemma}

We will assume from this point on that all curves in $\bndry H$ and meridian disks have been isotoped to have minimal intersection.

\begin{definition}
    A set $\mathbf D \subset \cal D$ is a \textbf{system of meridian disks} for $H$ if all the disks in $\mathbf D$ are pairwise disjoint, non-parallel and they cut $H$ up into  3-balls.
\end{definition}

If the genus of $H$ is $g$, then a minimal system of meridian disk is a system of $g$ meridian disks that cut $H$ up into a single ball.

\begin{definition}
    If $S$ is a punctured sphere and $\gamma\subset S$ is a properly embedded arc that is not boundary parallel, then $\gamma$ is a \textbf{wave} if both ends are in the one boundary component of $S$.
\end{definition}

Let $H$ be a handlebody, $\tc \subset \bndry H$ be a set of essential pairwise disjoint simple closed curves and $\mathbf D$ be a minimal system of disks for $H$.  From above we are assuming that the intersection between $\tc$ and $\mathbf D$ is minimal.  Let $S$ be the $2g$ punctured sphere produced when $\bndry H$ is cut along the boundary curves of the disks in $\mathbf D$.  Then $\tc_S = S \cup \tc$ is a set of properly embedded arcs in $S$.

\begin{definition}\label{defn: waveless system}
$\mathbf D$ is called \textbf{waveless} system of meridian disks with respect to $\tc$ if there are no waves in $\tc_S$.
\end{definition}

In \cite{C&R1} it is shown that given a set of essential pairwise disjoint simple closed curves in the boundary of a handlebody then a waveless system of disks can always be found.

\begin{definition}
    If $H$ is a handlebody and $\tc$ is a set of essential pairwise disjoint simple closed curves in $\bndry H$, then  $\tc$ meets the $\mathbf n$ \textbf{disk-condition} in $H$ if $|D|\geq n$ for every meridian disk $D$.
\end{definition}

See \cite{C&R1} for a sufficient and necessary condition for $\tc$ to satisfy the $n$ disk-condition in $H$.  It is also shown in \cite{C&R1} that there is an algorithm to decide if a set of curves in $\bndry H$ meets the $n$ disk-condition.

Next we are going to give a description of the construction of 3-manifolds that meet the `disk-condition'.  Let $H_1$, $H_2$ and $H_3$ be three handlebodies.  Let $S_{i,j}$, for $i\not= j$ be a sub-surface of $\bndry H_i$ such that:
\begin{enumerate}
    \item
        $\bndry S_{i,j} \not= \emptyset$.
    \item
        The induced map of $\pii(S_{i,j})$ into $\pii(H_i)$ is injective.
    \item
        For $j\not= k$, $S_{i,j}\cup S_{i,k} = \bndry H_i$,
    \item
        $\tc_i  = S_{i,j}\cap S_{i,k} = \bndry S_{i,j}$ is a set of disjoint essential simple closed curves that meet the $n_i$ disk-condition in $H$,
    \item
        and  $S_{i,j}\subset \bndry H_i$ is homeomorphic to $S_{j,i} \subset \bndry H_j$.
\end{enumerate}

Note that $S_{i,j}$ need not be connected. Now that we have the boundary of each handlebody cut up into essential faces  we want to glue them together by homeomorphisms, $\Psi_{i,j}: S_{i,j}\to S_{j,i}$, that agree along $\tc_i$'s.  The result is a closed 3-manifold $M$, for which the image of each handlebody is embedded.

\begin{figure}[h]
   \begin{center}
       \includegraphics[width=6cm]{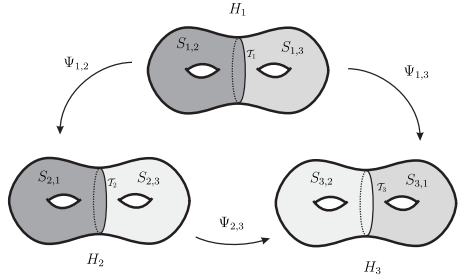}
    \end{center}
   \caption{Homeomorphisms between boundaries of handlebodies.} \label{fig:handlebodies}
\end{figure}

\begin{definition}
    If $M$ is a manifold constructed from three handlebodies as above, such that $\tc_i$ meets the $n_i$ disk-condition in $H_i$ and

    \begin{equation}\label{equation: disk-condition}
         \sum_{i=1,2,3} \frac{1}{n_i} \leq \frac{1}{2}
    \end{equation}

    then $M$ meets the $\mathbf{(n_1,n_2,n_3)}$ \textbf{disk-condition}. If we are not talking about a specific $(n_1,n_2,n_3)$, the manifold is said to meet just the \textbf{disk-condition}.
\end{definition}

Note that if the three handlebodies are solid tori then the resulting manifold is a non-Haken Seifert fibered manifold with infinite fundamental group. In \cite{C&R1} it is shown that the disk condition implies that the induced map on each $\pii(H_i)$ into $\pii(M)$ is injective and thus $\pii(M)$ is infinite. See \cite{C&R1} for further and possibly more enlightening discussion of this construction, in particular in relation to non-Haken Seifert fibered manifolds with infinite fundamental group and Haken  manifolds.  Also in \cite{C&R1} and \cite{Cof1} there are  a number of motivating examples. These included a construction for non-Haken non-Seifert fibered manifolds meeting the disk-condition. These are produced by performing Dehn Surgery along infinite families of pretzel knots and 2-bridge knots in $\s[3]$.

The following definition of CAT(k) metric spaces comes from the book \cite{b&h} by M.\,Bridson and A.\,Haefliger.  However, these spaces were first introduced  by Alexandrov in \cite{al1}. For the purposes of this paper we are only concerned with CAT(0) spaces, so the definition has been restricted to this one case.

If $X$ is a \textbf{length space}, it is a metric space such that any two points are joined by a unique geodesic and the distance between the points is the given by its length. Let  $p_1,p_2,p_3\in  X$ are three points and $[p_i,p_j]$ for $i\not= j$ be geodesic segments joining them. Then  $\triangle = [p_1,p_2]\cup [p_2,p_3] \cup [p_1,p_3]$ is a geodesic triangle.    Let $\bar p_1, \bar p_2, \bar p_3 \in \E[2]$ be three points such that $d(p_i,p_j) = d(\bar p_i, \bar p_j)$.   Then let the \textbf{comparison triangle} to $\triangle$ be   $\overline \triangle$  the geodesic triangle in $\E[2]$ with vertices $\bar p_1, \bar p_2, \bar p_3$ and the \textbf{comparison point} to $x\in [p_i,p_j]$ be the point  $\bar x \in [\bar p_i , \bar p_j]$ such that $d(p_i,x) = d(\bar p_i,\bar x)$.   $\triangle$ satisfies the CAT(0) \textbf{inequality} if for all points $x,y \in \triangle$ and comparison points $\bar x , \bar y \in \overline \triangle$, $d(x,y)\leq d(\bar x, \bar y)$

\begin{definition}
If the CAT(0) inequality holds for every geodesic triangle in a metric space $X$ then it is said to be CAT(0).
\end{definition}

A space is said to be \textbf{non-positively curved} if it is locally CAT(0).  That is, if every point $x\in X$ there is an $r > 0$ such that the ball of radius $r$ about $x$ is CAT(0). In \cite{b&h} Bridson and Haefliger show that a metric space being simply connected and non-positively curved is a sufficient condition for it to be CAT(0). All CAT(0) metric spaces are contractible, geodesic metric spaces.  Other properties of CAT(0) spaces required for this paper will be introduced as required.

\section{Structure in the universal cover}

\begin{figure}[h]
  \begin{center}
      \includegraphics[width=7cm]{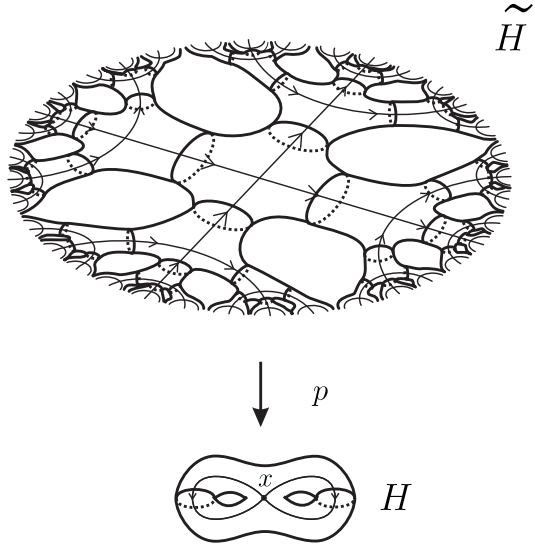}
   \end{center}
  \caption{Universal cover of a handlebody.}
  \label{fig:Model for the universal cover of a handlebody.}
\end{figure}

First we need to look at the universal cover of handlebodies.  Let $H$ be a handlebody with a set of curves $\tc\subset \bndry H$ that meet the $n$ disk-condition.  Let $\tld H$ be the universal cover of $H$ and $q:\tld H \to H$ be the covering projection.  Then, as can be seen in figure \ref{fig:Model for the universal cover of a handlebody.}, $H$ is P.L. 3-manifold with boundary, $\tld H$ is a non-compact P.L. simply connected 3-manifold with boundary and $int(\tld H) \homeo \R[3]$.  As each component of $\tc$ is essential in $H$, $q^{-1}(\tc)$ is a set of simple pairwise disjoint non-compact proper curves in $\bndry \tld H$.  Let $F\subset \bndry H$ be a face  produced when $\bndry H$ is cut along $\tc$.  By the disk-condition we know that $F$ is essential in $H$ and thus each component of $q^{-1}(F)$ is the universal cover of $F$.  Let $\tld F \subset \bndry \tld H$ be a component of $q^{-1}(F)$.  Then $\tld F$ is a simply connected non-compact P.L. 2-manifold, thus $int(\tld F) \homeo \R[2]$.

Note that  the universal covers of these objects has not been described by properly embedding them in compact manifolds. Even though this is a more concrete description it has not been done to avoid the complications brought about by ensuring that the compactification is `nice', as  tackled by Simon in \cite{si}. The following lemmas are required later in the paper when defining a process to glue the lifts of the handlebodies together to produce a sequence of expanding open 3-cells.

Let $M_1$ and $M_2$ be 3-manifolds such that $int(M_i)\homeo \R[3]$ and $\bndry M_i \homeo \R[2]$ and $f:\bndry M_1 \to \bndry M_2$ is a homeomorphism.  By R.\,Kirby \cite{Kr1} $f$ is a stable homeomorphism and by Brown and Gluck \cite{Br&GL1} all stable homeomorphism of $\R[n]$ are isotopic to the identity.  Then as shown by R.\,Lickorish in  \cite{lic}, the two manifolds resulting from gluings by isotopic homeomorphism, are homeomorphic.  This gives us the following lemma.

\begin{lemma} \label{lemma: interior of the union of two balls is an open ball}
    Let $M_i$'s and $f$ be described as above. If  $M_1 \cup_{ f} M_2 \homeo \R[3]$, then for any homeomorphism $f':\bndry M_1 \to \bndry M_2$, $M_1 \cup_{ f'} M_2 \homeo \R[3]$.
\end{lemma}

This means that whether the manifold resulting from the gluing is homeomorphic to $\R[3]$, depends only on the structure of the initial manifolds and not the homeomorphism between their boundaries that gives the gluing.  This leads us to the next lemma.

\begin{lemma}\label{lemma: condition for product structure}
    If the $M_i$'s are 3-manifolds as described above such that $M_1 \homeo \R[2]\cross [0,1)$ and $\bndry M_2$ has a collar neighbourhood in $M_2$, then for any homeomorphism $f:\bndry M_1 \to \bndry  M_2$, $M_1\cup_{f} M_2 \homeo \R[3]$.
\end{lemma}

\begin{proof}
Let $M= M_1\cup_{f} M_2$ and $U\subset M_2$ is the  collar neighbourhood of $\bndry M_2$. In an abuse of notation we will refer to the image of the $M_i$ in $M$ as $M_i$ and from the previous lemma we can assume that $f$ is the obvious identity.    It is sufficient to show that there is a homotopy of embeddings $r:M\cross I \to M$ such that $r(M,0) = M$ and $r(M,1)=int(M_2)$.  As $U$ is the  collar neighbourhood of $\bndry M_2$, then $U = \bndry M_2 \cross [0,1]$.    As $M_1 \homeo \bndry M_1 \cross [0,1)$ there is an obvious homotopy $r':(M_1\cup U)\cross I \to (M_1\cup U) $ where $r' (M_1\cup U,0)=M_1\cup U$, $r' (M_1\cup U,0)= int(U)$ and $r'$ restricted to $\bndry M_2\cross\{1\} $ is the identity.  Then $r$ is the identity on $\overline{M_2 -  U}$ and equivalent to $r'$ on $M_1\cup U$.
\end{proof}

Once again let $H$ be a genus $g$ handlebody, $\tld H$ its universal cover and $p:\tld H \to H$ is the covering projection. $\tc$ is a set of curves in $\bndry H$ meeting the $n$ disk-condition and $F$ is a face of $\bndry H$ when cut along $\tc$.

\begin{lemma}\label{lemma: product structure on lift of handlebody}
Let $\tld F\subset \bndry \tld H$ be a component of $p^{-1}(F)$, then $int(\tld H) \cup int(\tld F) \homeo \R[2] \cross [0,1)$.
\end{lemma}

\begin{remark}
This lemma seems reasonably obvious from the picture, however as it is a important step it is worth giving a formal proof.
\end{remark}

\begin{proof}
A necessary and sufficient condition for $N = int(\tld H) \cup int(\tld F) \homeo \R[2] \cross [0,1)$ is that any compact set in $N$ can be contained in a closed 3-cell. Therefore we can show that $N$ has a product structure by showing that for any compact set $C\subset N'$ there is a submanifold that contains $C$ and which admits has a product structure.

Let $p:\tld H \to H$ be the covering projection.  Let $\mathbf D$ be a minimal waveless system of meridian disks (see definition \ref{defn: waveless system}).  Let $\Gamma$ be the dual graph to $\mathbf D$, as shown in figure \ref{fig:Model for the universal cover of a handlebody.}. Therefore $\Gamma$ is a wedge of $g$ circles with a single vertex $x$. Note that each loop of $\Gamma$ corresponds to generator for representation of $\pii(H)$ and that $x$ is disjoint to $\mathbf D$. If $n(\Gamma)$ is the regular neighbourhood of $\Gamma$ in  $H$, then $H \homeo  n(\Gamma)$. Let $\tld \Gamma = p^{-1}(\Gamma)$, then $\tld \Gamma$ is an infinite tree and each vertex has order $2g$.  As  $\mathbf D$ is a set of $g$ meridian disks that cut $H$ up into a single 3-ball, then $\tld{\mathbf D} = p^{-1}(\mathbf D)$ is a set of pairwise disjoint properly embedded disks  that cut $\tld H$ up into fundamental domains, see figure \ref{fig:Model for the universal cover of a handlebody.}.  There is a single vertex of $\tld \Gamma$ contained in each fundamental domain.

Let $C$ be a compact set in $N$.  Then there is a simply connected  manifold $N'\subset \tld H$ made up of a finite number of fundamental domains of $\tld H$ such that $N'$ contains $C$.  Let $T\subset \tld \Gamma$ be the minimal connected tree such that $T$ contains all the vertices contained in $N'$. Therefore $T$ is a compact tree and $n(T) \subset int(H)$, its regular neighbourhood, is a closed three cell. Let $B$ be a fundamental domain in $N'$ such that $B\cap \tld F \ne \emptyset$. As $F$ is $\pii$-injective in $H$ and $\mathbf D$ is waveless, $B \cap \tld F$ is a disk. Let $\{ D_i\}$ be the maximal subset of $\tld{\mathbf D} \cap B$ such that $D_i \cap \tld F \ne \emptyset$.  $B$ can be fibered by $(B \cap \tld F) \cross [0,1]$ such that each $D_i$ is vertical.  Thus all the fundamental domains of $N'$ that have nonempty intersection with $\tld F$ can be similarly fibered so the fiberings agree on $\tld{\mathbf D}$ in their intersections.  This fibering can then be easily extended to the fundamental domains of $N'$ that have empty intersection with $\tld F$. We now have the product structure $(N'\cap \tld F) \cross [0,1]$ on $N'$  and thus the product structure $int(N'\cap \tld F) \cross [0,1)$ on $int (N') \cup int ( \tld F \cap N')$.  Thus $n(T)$ can be expanded in $int (N') \cup int ( \tld F \cap N')$ to contain $C$.
\end{proof}

Let $M$ be a 3-manifold which can be constructed by a gluing of three handlebodies, $H_i$'s, that meets the disk-condition and $\tc = \bigcap H_i$ is the set of triple curves.  Also let $X$ be the 2-complex $\bigcup \bndry H_i$. Note that this is not the usual use of 2-complex, as $X$ is constructed by gluing surfaces along their boundaries. These surfaces can easily be cut up into cells. Let $\tld M$ be the universal cover with $p:\tld M \to M$ the covering projection.  As the induced map of $\pii(H_i)$ into $\pii(M)$ is injective and  $H_i$ is embedded in $M$, each component of $p^{-1}(H_i)$ is the universal cover of $H_i$. As each face of $X$ is also $\pii$ - injective, $\tld X = p^{-1}(X)$ is a 2-complex whose faces are non-compact `ideal' polyhedron.  Note that if a face $X$ is not an annulus then its lifts in  $\tld X$ will have  infinite order as an ideal polygon.

\section{Dual 2-complex}

Next we wish to construct the dual 2-complex to $\tld X$.  There are two ways of viewing this object.  One is actually embedded in $\tld M$. The second is an abstract view where it exists by itself.  For the majority of this proof it does not matter which view we use, so most of the time no distinction will be drawn. However, each view does have its advantage at some point in this argument.

Let $\mathcal{C}$ be the dual 2-complex  to $\tld X$. The embedded view of $\mathcal C$ is that it has a vertex at the center of each lift of a handlebody.  There is an edge dual to each lift of a face of $X$ in $\tld M$ and a face dual to each lift of a curve in $\tc$.  As each curve in $\tc$ intersects three lifts of handlebodies,  each face of  $\mathcal C$ is a triangle.  The abstract view of $\mathcal C$ is that each  vertex  corresponds to a lift of a handlebody.  An edge between two vertices corresponds to the gluing, along the lift of a face of $X$,  between the two corresponding universal covers of handlebodies. Finally a face of $\mathcal{C}$ corresponds to the gluing of the lift of three handlebodies along the lift of a triple curve.  Clearly the resulting 2-complex in each case is identical. As all the faces in $\mathcal{C}$ correspond to the lifts of triple curves thus they are all triangular.

\begin{figure}[h]
  \begin{center}
      \includegraphics[width=5cm]{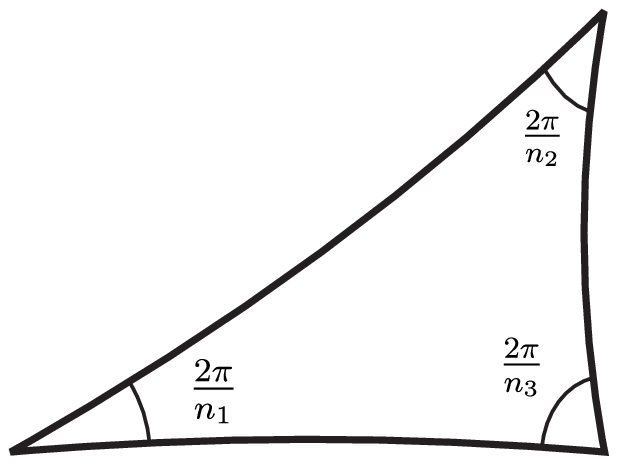}
   \end{center}
  \caption{A face of $\cal C$.}
  \label{fig: face of C.}
\end{figure}

We can put a metric on $\cal C$ by assuming all the faces are geodesics triangles, where the corner angles correspond to the disk-condition satisfied by the handlebody corresponding to the vertex.  If the vertex corresponds with a lift of the handlebody $H_i$ such that $\tc$ meets the $n_i$ disk-condition in $H_i$, then the internal angle of the face at the vertex is $2\pi/n_i$, as shown in figure \ref{fig: face of C.}. As each face has a vertex associated with each of the handlebodies, the sum of the internal angles of each face is always at most $\pi$.  Thus all the faces are similar and are either Euclidean or Hyperbolic triangles.  Note if the $M$ meets the $(6,6,6)$, $(4,8,8)$ or $(4,6,12)$ disk-condition then all the faces are  Euclidean geodesic triangles.  We will always assume that the shortest edge of each face has length 1. By Bridson and Haefliger \cite{b&h},  as the number of isometry classes of simplices of $\cal C$ is finite, $\cal C$ is a complete length space, meaning that the distance between any two points is the length of the shortest path joining them. If the genus of a handlebody is at least two then there are an infinite number of lifts of its faces in the boundary of a lift of the handlebody, thus the corresponding vertex in $\cal C$ is going to have infinite order.  Similarly if a face is not an annulus, all the edges in $\cal C$ corresponding to lifts of it will have infinite order. Let $\mathcal{C}^1$ be the 1-skeleton of $\mathcal{C}$ and $\mathcal C^0$ be the set of vertices of $\mathcal C$.

The simplest examples for $\cal C$ is when $M$ is a Seifert fibered space. As each handlebody is a solid torus, a vertex that is at the center of a lift of the solid torus $H_i$ has order $n_i$ and as each face of $X$ is an annulus every edge has order two.  Thus $\cal C$ is a tessellation of either $\E[2]$ or $\h[2]$ by isomorphic geodesic triangles, depending on the disk-condition.

Only the base disk-conditions ($(6,6,6)$, $(4,8,8)$ and $(4,6,12)$) need be considered.  The reason for this is that the disk-condition is used by the non-positively curved structure it induces and thus the minimal cases are somehow the worst. Also if $M$ satisfies the $(n_1,n_2,n_3)$ disk-condition, such that $\sum \frac{1}{n_i} < \frac{1}{2}$, then $M$ also meets at least one of the base disk-conditions.  This means that we can assume the faces of $\cal C$ are Euclidean geodesic triangles.


$\cal C$ seems fairly unfriendly as it can have infinite degree on its vertices and edges, however it turns out not too bad as it is CAT(0).  From Bridson and Haefliger \cite{b&h} a length space being simply connected and non-positively curved is  sufficient to show that $\cal C$ is CAT(0).

\begin{lemma}
If $M$ is a manifold that meets the disk-condition, $\tld M$ is its universal cover and $\cal C$ is the dual 2-complex to $\tld X$, as described above, then $\cal C$ is simply connected.
\end{lemma}

\begin{proof}
The idea of this proof is to show that there is a retraction of $\tld M$ onto $\cal C$, therefore as $\tld M$ is simply connected then $\cal C$ is simply connected.  Let $\tld \tc = p^{-1}(\tc)$ and $\tld X =  p^{-1}(X)$.  Through out this proof we will refer to the embedded model of $\cal C$.  We will also assume that it is in general position with respect to $\tld X$.

First we want to show there is a retraction of $\tld \tc$  into $\cal C$. Each component of $\tld \tc$ is an open arc that intersects $\cal C$ exactly once in the interior of a 2-simplex.  Moreover each face intersect exactly one component of $\tld \tc$.  This means  that for any face $F$ of $\tld X$ that $F\cap \cal C$ is a tree that cuts $F$ up into 4-gons with a single vertex removed, as shown in figure \ref{fig: retraction.}.  There is a retraction $r_{\tld \tc} : \tld \tc \to \cal C$ where each component $\alpha \subset \tld \tc $ is sent to the point $\alpha \cap \cal C$.  That is $r_{\tld \tc}(\tld \tc) = \tld \tc \cap \cal C$

Next we want to show there is a retraction of the lifts of the faces of $X$ to $\cal C$.  Let $A$ be a face of $\cal C$ and $A'$ a face adjacent  to $A$, such that $a=A\cap \tld \tc $ and $a' =A'\cap \tld \tc $ can be joined by a simple arc in $\tld M$ otherwise disjoint from $\cal C$.  Let $\gamma$ and $\gamma'$ be the open arcs from $\tld \tc $ such that $\gamma\cap A = a$ and $\gamma'\cap A = a'$.  Let $\beta$ be the simple sub-arc of the graph $\tld X \cap \cal C$ that joins $a$ and $a'$. Thus $\beta$ is contained in a single face $F$ of $\tld X$ and cuts  $F'$ from $F$ which is a quadrilateral, with a vertex removed.  Thus the retraction of $\gamma$ to $a$ and $\gamma'$ to $a'$ can be extended to a retraction of $F'$ to $\beta$.  Thus there is a retraction of $F$ to the tree $\cal C \cap F$.  We can then extend this to a retraction $ r_{\tld X}:\tld X \to \cal C$, where $r_{\tld X} (\tld X) = \tld X \cap C$, which is the dual graph to $\cal C$.

\begin{figure}[h]
  \begin{center}
      \includegraphics[width=6cm]{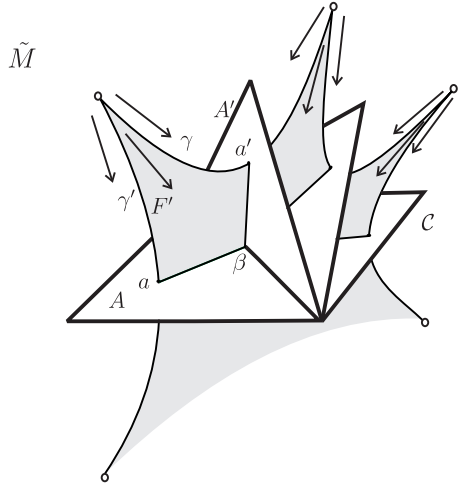}
   \end{center}
  \caption{Retraction of $\tld X$ into $\cal C$.}
  \label{fig: retraction.}
\end{figure}

As we get $\tld M$ from $\tld X$ by gluing in open 3-cells to $\tld X\cup \cal C$, the retraction $r_{\tld X}$ can be extended to a retraction $r:\tld M \to \cal C$. Thus $\cal C$ and $\tld M$ are homotopy equivalent and $\cal C$ is simply connected.
\end{proof}


\begin{definition}
    For $x\in \mathcal C$ let the \textbf{star} of $x$, denoted $st(x)$, be the union of all simplices containing $x$. Let the \textbf{open star}, denoted $st^o(x)$, be the union of $x$ and of the interiors of all the simplices adjacent to $x$. Then the \textbf{link} of $x$, denoted $lk(x)$,  is the union of all simplices in $st(x)$ that do not contain $x$.
\end{definition}

\begin{definition}
For simplices  $A,B\subset \cal C$, if $A\subset B$ then $B$ is said to be \textbf{adjacent} to $A$.
\end{definition}

\begin{definition}
 If $\cal D$ is a simplicial complex and $e\subset \cal D$ is an edge, then $e$ is said to be a \textbf{free edge} if it is adjacent to exactly one face.
\end{definition}

Note that $\cal C$ does not have any free edges, and thus for $x\subset \cal C$ some vertex, $lk(x)$ is the union of all the free edges of $st(x)$.

\begin{lemma}\label{lemma: C has non+'ve curvature}
$\cal C$ is non-positively curved.
\end{lemma}

As each face is a non-positively curved geodesic triangle, if  $\cal C$ is going to have positive curvature then it is going to be at a vertex.  From Bridson and Haefliger \cite{b&h} showing that every non-trivial path in the link of a vertex has length at least $2\pi$  is a sufficient condition for $\cal C$ to be non-positively curved.

\begin{proof}
For $v\in \cal C^0$, we can think of $lk(v)$ as a graph in $S^2$,  where each edge of $lk(v)$ represents a face in $\mathcal{C}$ and a vertex in $lk(v)$ represents an edge in $\mathcal{C}$ .  If $v$ corresponds to a lift of a handlebody that meets the $n$ disk-condition then we can give each edge in $lk(v)$ a length of $\frac{2\pi}{n}$.  This means that the length of an edge in $lk(v)$ is the same as the angle between the other two edges in $\cal C$ of the triangle represented by the edge.  If $\alpha$ is a loop in $lk(v)$, then let $|\alpha|$ be the length of $\alpha$. Let $\tld v$ be the lift of a handlebody in $\tld M$ corresponding to $v$. Using the embedded view of $\cal C$, a loop $\alpha'$ in $\Gamma = \cal C \cap \bndry \tld v \homeo lk(v)$,  that is essential in $\Gamma$, and thus  $\bndry \tld H$, will project down to a curve that is the boundary of a singular meridian disk i $H$.  Thus by the disk condition $\alpha'$ will intersect $\cal C^1$ at least $n$ times. From the obvious retraction of $st(v)$ onto $\cal C \cap \tld v$, for any non-trivial loop $\alpha \subset lk(v)$, $|\alpha |\geq 2\pi$.   From Bridson and Haefliger \cite{b&h} we know this is a sufficient condition for $\mathcal{C}$ to be non-positively curved.
\end{proof}

A direct corollary of the previous two lemmas is:

\begin{corollary}
If $\cal C$ is  the 2-complex described above, then it is CAT(0).
\end{corollary}

\section{Expanding sequence in $\cal C$}

As $\cal C$ is CAT(0) any two  points in $\cal C$ can be joined by a unique geodesic.   Also as $\cal C$ has no free edges, that is all edges are adjacent to at least two faces, any geodesic arc can be extended. Thus we can define a metric ball in $\cal C$ around some vertex $ v \subset \cal C ^0$.


\begin{definition}\label{defn: metric balls et al}
    Let a \textbf{metric ball} of radius $r$ centered at $v$ be  $\mathbf B_{r,v} =\{x\in \cal C : d(x,v)\leq r\}$ and a \textbf{level sphere} be $\mathbf S_{r,v} =\{x\in \cal C : d(x,v)= r\}$.   A level sphere is said to be \textbf{regular} if it is transverse to the simplicial structure of $\cal C$, otherwise it is said to be \textbf{critical}. Also for $\mathbf S_{r,v}$ critical, let $V_{r,v}=\mathbf S_{r,v} \cap \cal C^0$.
\end{definition}

As $\cal C$ is CAT$(0)$, $\mathbf B_{r,v}$ is convex and contractible. Also note that if a level sphere is regular then it is disjoint from $\cal C^0$.  A level sphere is a connected graph, that is separating in $\cal C$ and whose edges are circular arcs.  Therefore as the edges of $\cal C$ are geodesic, for $r>0$, $\mathbf B_{r,v}$ is never simplicial.  If there was a sequence of simplicial metric balls, then they would correspond to  a sequence of expanding  submanifolds in $\tld M$ and thus all that would remain to prove theorem \ref{theorem: universal cover is $R^3$} is that the interior of these manifolds is homeomorphic to $\R[3]$. Unfortunately this is not the case and further work is required to generate an expanding sequence of simplicial  subsets of $\cal C$.

As any geodesic arc can be extended, we can produce an infinite sequence of expanding metric balls in $\cal C$.  We do this by showing that the radii of critical spheres around a vertex are discrete. This means that if we take the infinite sequence of expanding critical spheres about a vertex, that the union of the sequence is $\cal C$.

\begin{lemma}\label{lemma: critical shperes don't accumulate}
    Critical level spheres do not accumulate at a finite radius $r$.
\end{lemma}

\begin{remark}
A  version of this lemma for CAT$(0)$ polyhedral complexes is proved by Bridson and Haefliger  in \cite{b&h},  pg 122.
\end{remark}

\begin{proof}
The idea of this proof is to show that the distances between any two vertices in $\cal C$ form a discrete set in $\R[]$.  That is any range in $\R[]$ contains only finitely many possible distances.  Let $\triangle$ be the Euclidean geodesic triangle isometric to each face of  $\cal C$ and $\mathbf T$ be the tessellation of $\E[2] $ by $\triangle$.  Therefore any disk of finite radius contains a finite number of vertices.  This means that if $D$ is the set of all possible distances between vertices in $\mathbf T$, then there are only finitely many elements of $D$ smaller than any finite $r>0$.  In other words $D$ is discrete in $\R[]$. Also note that as the shortest edge length of $\triangle$ is always 1, then for any $d\in D$ $d\geq 1$.

Let $\gamma$ be a geodesic  in $\cal C$, between two vertices and $l(\gamma)$ be its length.  Note that $l(\gamma)$ is finite.  If $\gamma$ is disjoint to $\cal C^0$ other than at its ends, then it is isometric to a geodesic in $\E[2]$ between vertices of $\mathbf T$.  This means that the possible lengths of such a geodesic are discrete.  If $\gamma$ intersect vertices of $\cal C$  other than at its ends, it must do so at most $l(\gamma)$ times, as the minimum distance between vertices is 1.   Each sub-geodesic of $\gamma$ which is disjoint from $\cal C^0$ other than at its ends is isometric to a geodesic between vertices of $\mathbf T$.  Therefore if $l(\gamma)\leq r$, $l(\gamma) = \sum_{i=1}^m d_i$ where $d_i\in D$ and $m\leq r$.   However, as $l(\gamma)\leq r$ each $d_i\leq r$.  There are only a finite set of such $d_i$'s.  Thus, as there is only a finite number of  finite sums from a finite set of numbers, the possible lengths of $\gamma$ are finite if its length is at most $r$.   This means that  the possible lengths of geodesics between vertices in $\cal C$ are a discrete set.
\end{proof}

\begin{definition}\label{defn: neighbouring vertices}
    Two vertices $v,v'\in\cal C$ are said to be \textbf{neighbours} if they are both adjacent to the same edge.
\end{definition}

\begin{lemma}\label{lemma: st(v) is convex}
    For a vertex $v\subset \cal C$, $st(v)$ is convex.
\end{lemma}

\begin{proof}
Clearly $st(v)$ is connected as any two points in it can be joined by a path through $v$.  Note that  $lk(v)$ is the union of all the free edges of $st(v)$.  As the shortest edge length in $\cal C$ is always 1, $\mathbf S_{\frac{1}{2},v}$ is contained in $st(v)$.  Moreover, $lk(v)$ is homeomorphic to $\mathbf S_{\frac{1}{2},v}$. Therefore we know that $lk(v)$ is connected and separating.  Let $l$ be the length of the longest edges adjacent to $v$.  This means that if $M$ meets the $(6,6,6)$ disk-condition, then $l = 1$, if it meets the $(4,6,12)$ disk-condition, then $l = \sqrt{3}$ or $ 2$ and if it meets the $(4,8,8)$ disk-condition then $l= 1$ or $\sqrt{2}$.  Let $A= st(v)\cap \mathbf S_{l,v}\subset lk(v)$, that is, all the vertices in $st(v)$ distance $l$ from $v$.  If all the edges adjacent to $v$ are  the same length then $A$ contains all the vertices of $lk(v)$.  If the edges adjacent to $v$ are two different lengths then, as each handlebody in $H$ is embedded, if a vertex in $lk(v)$ is in $A$ then all its neighbouring vertices in $lk(v)$ are not and vice versa. Thus $A$ is separating in $\overline{B_{l,v}-st(v)}$.

If $st(v)$ is not convex it must contain two points, $x$ and $y$, such that $\alpha$, the image of the geodesic joining them, is not contained entirely in $st(v)$. As $lk(v)$ is separating, $\alpha$ must intersect it at least twice and an even number of times.  We will also assume that $\alpha$ is transverse to $\cal C^1$.  As $B_{l,v}$ is convex and $A$ is separating in $\overline{B_{l,v}-st(v)}$ there is $\alpha ' $, a connected sub-arc of $\alpha$, such that $\alpha'\cap st(v) = \{x',y'\} \subset lk(v)$ where $x'$ and $y'$ are points contained in a single component of $lk(v) - A$.   This means they are either contained in the same edge of $lk(v)$ or in two adjacent edges of $lk(v)$.

Both $x'$ and $y'$ cannot be contained in the same  edge, as each edge is geodesic this would contradict the property of CAT$(0)$ spaces that any two points are connected by a unique geodesic.  Thus $x'$ and $y'$ must be contained in adjacent edges of $lk(v)$, $E_{x'}$ and $E_{y'}$ respectively.  Therefore $z=E_{x'}\cap E_{y'}\subset lk(v)$ is a vertex.  Let $E_z$ be the edge between $z$ and $v$.  As both $E_{x'}$ and $E_{y'}$ are free edges of $st(v)$ they are contained in the boundary of a single face of $st(v)$, $F_{x'}$ and $F_{y'}$ respectively, where $F_{x'}\cap F_{y'} = E_{z}$.  Thus $F_{x'}\cup F_{y'}$ is convex and contains a geodesic between $x'$ and $y'$.  This once again contradicts the CAT$(0)$ structure of $\cal C$. Thus $st(v)$ is convex.
\end{proof}

\begin{definition}
    For any set $C\subset \cal C$ and a point $v\not\in C$, let $d_C(v)= inf\{d(y,v): y\in C\}$.
\end{definition}

As $\cal C$ is a simplicial 2-complex and its faces are Euclidean triangles there is a simplified definition for the angles between geodesics. Let $x,y,z \subset \cal C$ be points in  $\cal C$ and  $[x,y]$ be the geodesic arc joining the points $x$ and $y$. Let the \textbf{angle} between $[x,y]$ and $[x,z]$ be denoted by $\angle_x(y,z)$.  If $x$ is not a vertex of $\cal C$ the  $\angle_x(y,z)$ is simply the Euclidean angle between $[x,y]$ and $[x,z]$.  If $x$ is a vertex of $\cal C$ then $\angle_x(y,z)$ is the length of the shortest path in $lk(x)$ between $[x,y]\cap lk(x)$ and $[x,z]\cap lk(x)$, as defined in the proof of lemma \ref{lemma: C has non+'ve curvature}.

From Bridson and Haefliger \cite{b&h} we get the following lemma,

\begin{lemma}\emph{(Bridson and Haefliger)}\label{lemma: B&H convex complete sets}
Let $C\subset \cal C$ be convex complete and  $v \not\in \cal C$ any point, then:
\begin{enumerate}
\item $\{y\in C:d(y,v) = d_C(v)\}$ is a unique point $x$,
\item if $y\in C$ and if $y\not= x$ then $\angle_x(y,v) \geq \pi/2$, where $x$ is as in $1$.
\end{enumerate}
\end{lemma}

As faces of $\cal C$ are convex, this lemma tells us two things. First that the subset of a face of $\cal C$ closest to any vertex is a single point in its boundary.  Secondly that if the closest point to a vertex is contained in the interior of an edge, then the geodesic between the point and the vertex intersect is normal to the edge.  This may seem rather obvious, however it can help us classify the types of intersections between a metric ball and faces.

\begin{lemma}\label{lemma: at least one y_i is a vertex}
Let $F\subset \cal C$ be a face and $v$ a vertex not in $F$.  Let $E_i$, for $i=1,2,3$, be the three edges of $F$ and $y_i=\{z\in E_i:d(z,v)=d_{E_i}(v)\}$. Then at least one $y_i$ is a  vertex of $F$.
\end{lemma}

\begin{proof}
The idea to this proof is to show that by assuming the lemma is false, we can find geodesics that intersect twice, giving us a contradiction.

Assume that the lemma is not true and that $y_i\in int(E_i)$ for $i=1,2,3$. Note by lemma \ref{lemma: B&H convex complete sets} we know that each $y_i$ is a single point and by assumption they are all different. Let $y$ be the point in $F$ closest to $v$. As $v\not\subset int(F)$, then $y$ is one of the $y_i$'s.  Let it be $y_1$, thus $d(y_1,v) < d(y_i,v)$ for $i=2$ or $3$.  Let $\gamma_i$ be the geodesic between $y_i$ and $v$.  By  lemma \ref{lemma: B&H convex complete sets} we know that $\gamma_i$ is perpendicular to $E_i$. Let $F'$ be the other face adjacent to $E_1$ such that $\gamma_1\cap F'\not= \emptyset$.  Let $E_i'$ be the edge of $F'$ isometric to $E_i$, thus $E_1 = E_1'$.  Two vertices of $F'$ are $x_1$ and $x_2$, let the third be $x_3'$.  As $y_i\subset int(E_i)$ and $y_1$ is the point of $F$ closest to $v$, then $\mathbf B_{d(y_1,v),v} \cap F = y_1$.  Assume that $\gamma_i \cap E_1 = \emptyset$, for $i=2$ or $3$. Then, as any two points in $\cal C$ are connected by a unique geodesic arc, $B_{d(y_2,v),v} \cap F$ must  have two components, one of which is $y_i$ and the other has $y_1$ in its interior.  This contradicts convexity, thus both $\gamma_2$ and $\gamma_3$ intersect $int(E_1)$ and $int(F')$.

\begin{figure}[h]
  \begin{center}
      \includegraphics[width=6cm]{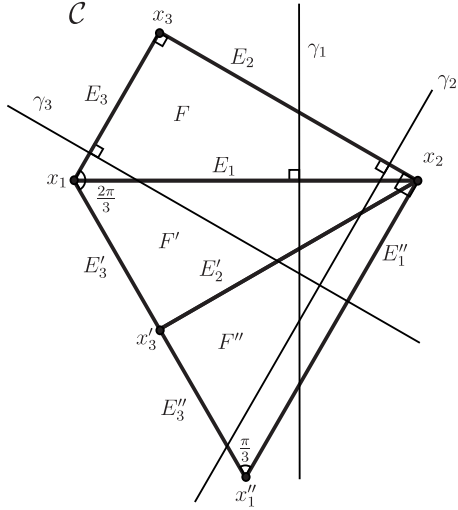}
   \end{center}
  \caption{Intersecting geodesics.}
  \label{fig: intersecting geods.}
\end{figure}

Extend $\gamma_1$ through $F$ until it intersects at least one of $E_2$ or $E_3$. Assume that it intersects $E_2$. As $\gamma_1$ is perpendicular to $E_1$ it also intersects $E_2'$. However if $\angle_{x_2}(x_1,x_3)\geq \pi/4$ then $\angle_{x_2}(x_3',x_3)\geq \pi/2$, meaning that  $\gamma_2$ intersects $E_3'$. Thus $\gamma_1\cap \gamma_2 \not= \emptyset$, which contradicts the unique geodesic property of CAT(0) spaces.  Therefore $\angle_{x_2}(x_1,x_3)= \pi/6$. This means that the length of $E_1$ is either 2 or $\sqrt{3}$.  If $E_1$ has length $\sqrt{3}$, then $\angle_{x_1}(x_2,x_3)= \pi/2$ and $\gamma_3\cap F$ is parallel to $E_1$ and thus would intersect both $\gamma_1$ and $\gamma_2$.  Therefore the length of $E_1$ must be $2$ and $\angle_{x_1}(x_2,x_3)= \pi/3$. Thus $\angle_{x_1}(x_3',x_3)= 2\pi/3$.  This means that $\gamma_3 \cap E_3'=\emptyset$ and thus $\gamma_3 \cap E_2' \not= \emptyset$.  This means that each $\gamma_i$ intersects $int(E_2')$.  Let $F''$ be the other face adjacent to $E_2'$ such that $\gamma_1\cap int(F'')\not= \emptyset$. Thus $\gamma_i\cap int(F'')\not= \emptyset$ for each $\gamma_i$.  Once again let $E_i''$ be the edge of $F''$ isometric to $E_i$ and let the vertex $E_3''\cap E_1'' = x_1''$.  Note that $E_3''\cap E_3'$ is a geodesic and thus $\angle_{x_1}(x_1'',x_3)= 2\pi/3$.  Therefore $E_1''\cap \gamma_3\not=\emptyset$.  Also $\angle_{x_2}(x_1'',x_3)= \pi/2$, thus  $(E_3'\cup E_3'')\cap \gamma_2\not=\emptyset$.  Thus $\gamma_2$ and $\gamma_3$ intersect in $int(F\cup F' \cup F'')$.  This takes care of all the cases, thus proving the lemma.
\end{proof}

A direct corollary of lemma \ref{lemma: at least one y_i is a vertex} is:

\begin{corollary}\label{corollary: S_rv intersect F can't have three components}
    If $F$ is a face of $\cal C$ and $v$ is a vertex, then for any $r> 0$, $\mathbf S_{r,v}\cap F$ has at most two components.
\end{corollary}

\begin{figure}[h]
  \begin{center}
      \includegraphics[width=3cm]{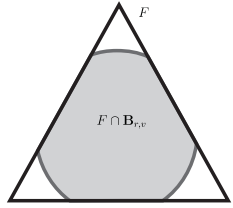}
   \end{center}
  \caption{$F \cap \mathbf B_{r,v}$ if $F \cap \mathbf S_{r,v}$ has three components.}
  \label{fig: S_rv_has_three_arcs_in_F.}
\end{figure}

\begin{proof}
Suppose the corollary is not true and that for some face $F\subset \cal C$, $F\cap S_{r,v}$ has three components. By convexity of $F$ and $\mathbf B_{r,v}$, $F \cap \mathbf B_{r,v}$ has to be connected and convex, thus by convexity we know that $F \cap \mathbf B_{r,v}$ looks like figure \ref{fig: S_rv_has_three_arcs_in_F.}, with $F \cap \mathbf S_{r,v}$ consisting of three properly embedded arcs, one running between each pair of edges of $F$.  This contradicts lemma \ref{lemma: at least one y_i is a vertex}, for it implies that the closest point to $v$ for each edge is in their interior.
\end{proof}

\begin{lemma}\label{lemma: S_r can't have parallel components in F}
    If $F\cap \mathbf S_{r,v}$  has two components then they cannot connect the same pair of edges.
\end{lemma}

\begin{proof}
Once again let $x_i$, for $i=1,2,3$, be the vertices of $F$ and $E_i$ be the edge between $x_i$ and $x_{i+1}$, where $E_3$ is the edge between $x_3$ and $x_1$. Also let $\{y_i\}=\{z\in E_i:d(z,v)=d_{E_i}(v)\}$.  Let $\gamma_i$ be the geodesic that passes through $y_i$ and $v$.  Note that the geodesic $\gamma_i$ is perpendicular to $E_i$.  If $F\cap \mathbf S_{r,v}$ connects the same pair of edges then at least two of the $y_i$'s are in $int(E_i)$.  Let $y_1\subset int(E_1)$  and  $y_2\subset int(E_2)$.  Using the same argument, as was used in the proof of lemma \ref{lemma: at least one y_i is a vertex}, $\gamma_2\cap E_1\not= \emptyset$. Also $d(y_1,v)<d(y_2,v)\leq d(y_3,v)$.  This means that $\alpha = \mathbf S_{d(y_2,v),v}\cap F$ is a circular arc with both ends in $E_1$ and $\alpha \cap E_2=y_2$.

By simple geometry if $M$ meets either the $(6,6,6)$ or $(4,8,8)$ disk-conditions, then the radius of $\alpha$ has to be less than 1. This  would imply that   $\gamma_2$ intersects a vertex $x\subset \mathbf B_{d(y_2,v),v}$ such that $d(y_2,x)<1$.  The only such vertex is $x_1$, however that would contradict $d(y_2,v)\leq d(y_3,v)$.

If $M$ meets the $(4,6,12)$ disk-condition and either $E_1$ or $E_2$ have length 1, then the radius of $\alpha$ has to be less than 1. This would once again mean that $x$ is $x_1$. Let $F'$ be the other  face adjacent to $E_1$ such that $\gamma_1\cap int(F')\not= \emptyset$ and let $x'$ be the vertex of $F'$ disjoint to $F$.  Then  $\gamma_2\cap F'\not= \emptyset$.   If $E_1$ has length $\sqrt{3}$, then $E_2$ must have length 2.  In this case then the radius of $\alpha$ must be less than $2/\sqrt{3}$,  once again this would mean that $x$ is $x_1$.  Therefore $E_1$ must have length 2 and $E_2$ has length $\sqrt{3}$.  In this case $\alpha$ must have radius less than $2$.  Therefore $x$ is either $x_1$ or $x'$.  Once again it cannot be $x_1$. If $x=x'$ then $d(x,x_1)=1$, once again contradicting that $d(y_2,v)\leq d(y_3,v)$.  Thus having looked at all the cases the lemma has been proved.
\end{proof}


\begin{lemma}\label{lemma: getting epsilon small enough}
    Let  $r>0$  where $\mathbf S_{r,v}$ is critical, then for  $0<\epsilon < r- \sqrt{r^2 - \frac{1}{4}}$  the level sphere $S_{r-t, v}$, where $t\in (0,\epsilon ]$, transversely intersects each edge between any vertex in $V_{r,v}$ (see definition \ref{defn: metric balls et al}) and any of its neighbouring vertices contained in $\mathbf B_{r,v}$.
\end{lemma}

\begin{figure}[h]
  \begin{center}
      \includegraphics[width=3cm]{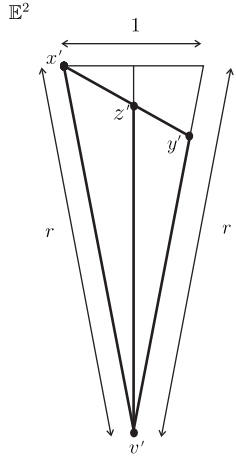}
   \end{center}
  \caption{Comparison triangle.}
  \label{fig: Comparison triangle.}
\end{figure}

\begin{proof}
The idea of this proof is to show that for any  vertex $x\in V_{r,v}$ and an edge $E\subset \mathbf B_{r,v}$ adjacent to $x$, that there is at least one point  $z \in E$ such that $d(z,v)\leq \sqrt{r^2 - \frac{1}{4}}$.  Let $y\subset \mathbf B_r$ be the vertex at the other end of $E$.  Let $\gamma_x$ and $\gamma_y$ be the geodesic between $v$ and $x$ or $y$ respectively. Then as $E$ is a geodesic, $\triangle = E \cup \gamma_x \cup \gamma_y$ is a geodesic triangle.

Let $x'$, $y'$ and $v'$ be points in $\E[2]$ such that $d(x',v') = d(x,v)=r$, $d(y',v') = d(y,v)\leq r$ and  $d(x',y') = d(x,y)$. For points $a$ and $b$ in $\E[2]$,  let $\overline{ab}$ be the image of the geodesic arc between $a$ and $b$.  Therefor $\triangle ' = \overline{x'v'} \cup \overline{y'v'}\cup \overline{x'y'}$ is a geodesic triangle in $\E[2]$.  We refer to $\triangle '$ as the comparison triangle to $\triangle$, see figure \ref{fig: Comparison triangle.}.  As $\cal C$ is CAT$(0)$ we know that for any point $z\subset E$ and $z'\subset \overline{x'y'}$, such that $d(z,x)=d(z',x')$, then $d(z,v)\leq d(z',v')$. Let $z$ be the point half way along $E$, that is $d(x,z)=\frac{1}{2}d(x,y)\geq\frac{1}{2}$, as the shortest edge length in $\cal C$ is always $1$ .  Let $z'\subset\overline{xy}$ be the corresponding point in $\triangle '$, that is $d(z',x')=d(z,x)$. $d(z',v')$ is maximal when $d(y',v')=r$ and $d(x',y')=1$,  therefore $d(z,v)\leq d(z',v')\leq \sqrt{r^2 - \frac{1}{4}}$.
\end{proof}

So far we can produce a sequence of expanding balls in $\cal C$.  However for $r>0$ these balls are not simplicial, thus there is no corresponding manifold in $\tld M$, for there is no obvious way to have part of a gluing along a face or around a triple curve.  Thus we need to produce a sequence of expanding simplicial sets.

\begin{definition}
Let the \textbf{simplicial ball} of radius $r$ centered at $v$, denoted $\mathbf B_{r,v}^s$ be the maximal simplicial subset of $\mathbf B_{r,v}$, the metric ball of radius $r$.
\end{definition}

If $M$ meets the $(4,8,8)$ disk-condition and $v \subset \cal C^0$ is a vertex such that it has adjacent edges of length 1 and $\sqrt{2}$.  Then, for $1\leq r <\sqrt{2}$, $\mathbf B_{r,v}^s$ is the set of all edges adjacent to $v$ with length 1.  This is not a convex set in $\cal C$. Therefore the name simplicial ball is a little misleading, as these objects are not necessarily convex.

\begin{lemma}\label{lemma: retraction of B_r onto B_r^s}
    There is a retraction of $\mathbf B_{r,v}$ onto $\mathbf B_{r,v}^s$.
\end{lemma}

\begin{proof}
The idea of this proof is to show that any simplex not entirely contained in $\mathbf B_{r,v}$  can be retracted into its boundary. Let $A\subset \cal C$ be a simplex such that $A\cap \mathbf B_{r,v}\not= \emptyset$ and $A\cap \mathbf B_{r,v}\not= A$.  Then as both $A$ and $\mathbf B_{r,v}$ are convex and $A$  is compact, then $A\cap \mathbf B_{r,v}$ must be convex, compact and thus connected.  Therefore if $A$ is an edge there are two possible cases.  The first is that $A\cap \mathbf B_{r,v}$ is a compact connected set in int$(A)$, the second is that a vertex of $A$ is contained in $A\cap \mathbf B_{r,v}$.

If $A$ is a face then the arcs $A\cap \mathbf S_{r,v}$ are circular arcs. By lemma \ref{corollary: S_rv intersect F can't have three components} we know that  $A\cap \mathbf S_{r,v}$ has at most two components. If $A\cap \mathbf S_{r,v}$ is a single arc running between different edges of $A$, then type 1 and type 2 in figure \ref{fig:Retraction of  B_r to B_r^s.} are the two possibilities for $A\cap \mathbf B_{r,v}$.  If $A\cap \mathbf S_{r,v}$ is a single arc with both ends in the one edge $E$, then by convexity we know that $A\cap \mathbf B_{r,v}$ must be the disk between the arc $A\cap \mathbf S_{r,v}$ and $E$, as shown in type 3 of  figure \ref{fig:Retraction of  B_r to B_r^s.}.  If $A\cap \mathbf S_{r,v}$ contains two arcs, by convexity neither can be homotopic to an edge of $A$ and by lemma \ref{lemma: S_r can't have parallel components in F} they cannot connect the same edges.  Therefore $A\cap \mathbf S_{r,v}$ must look like type 4 of figure \ref{fig:Retraction of  B_r to B_r^s.}.

\begin{figure}[h]

     \subfigure[type 1]{\hspace{2cm}
                \includegraphics[width=4cm]{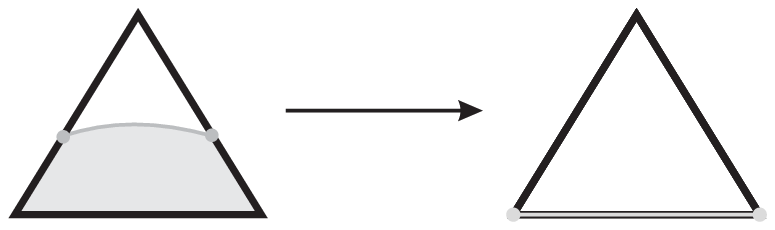}}

      \subfigure[type 2]{\hspace{2cm}
                \includegraphics[width=4cm]{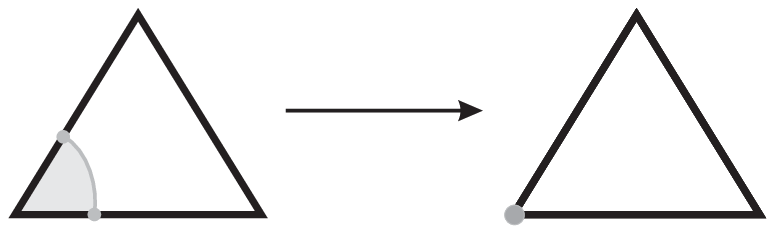}}

      \subfigure[type 3]{\hspace{2cm}
                \includegraphics[width=4cm]{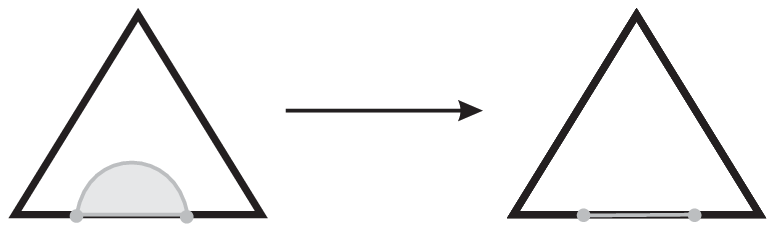}}

      \subfigure[type 4]{\hspace{2cm}
                \includegraphics[width=9.3cm]{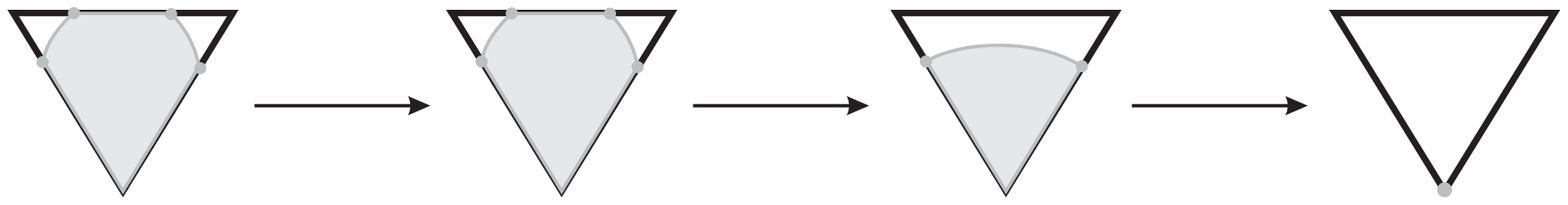}}

  \caption{Retraction of $\mathbf B_r$ to $\mathbf B_r^s$.}
  \label{fig:Retraction of  B_r to B_r^s.}
\end{figure}

Let $E\subset \cal C$ be an edge such that  $E\cap \mathbf B_{r,v} \subset \mathrm{int}(E)$.  Then if $F\subset \cal C$ is a face adjacent to $E$, $F\cap \mathbf S_{r,v}$ either has one arc parallel to $E$ or has two arcs.  Therefore $F\cap \mathbf B_{r,v}$ is either of type 3 or type 4.  However $E$ cannot have more than one adjacent face whose intersection with $\mathbf B_{r,v}$ is of type 4, for if it did then any point in $E\cap \mathbf B_{r,v}$ would have more then one geodesic path to $v$.  However all the adjacent faces of $E$ cannot have intersections of type 3 with $\mathbf B_{r,v}$ for the circular arcs of $\mathbf S_{r,v}$ are centered at a vertex in $\cal C$.   This means that $E$ must have one adjacent face that has an intersection of type 4 with $\mathbf B_{r,v}$ and all its other adjacent faces must have intersection of type 3 with $\mathbf B_{r,v}$.

Therefore for any face $F$ adjacent to $E$ that has an intersection of type 3 with $\mathbf B_{r,v}$ there is a retraction of $\mathbf B_{r,v}\cap F$ into $E$.  We then do this for all such faces adjacent to $E$.  Now we can retract the resulting set, $B$, to remove any intersection with $E$.  This means that  the one face adjacent to $E$ that has an intersection of type 4 with $\mathbf B_{r,v}$ now has an intersection of type 2 with $B$.  We can then repeat this process to remove all intersections of type 3 and type 4 from $B$.

This means that for any face  $F\subset \cal C$ not entirely contained in $B$, that the intersection  between them is either of type 1 or type 2. Clearly in either case $B$ can be retracted so that  it only intersects $F$ in its boundary.  Now the only simplices that are not entirely in $B$  are edges with one vertex contained in $B$.  Therefore we can retract $B$ so that its intersection with such edges is just the vertex.  Therefore $B = \mathbf B_{r,v}^s$ is now simplicial.
\end{proof}

Note that if a vertex is contained in $\mathbf B_{r,v}$ then it is in $\mathbf B_{r,v}^s$ and if two neighbouring vertices are contained in $\mathbf B_{r,v}$, then the edge joining them is contained in $\mathbf B_{r,v}^s$.  The next corollary is a direct result of lemmas \ref{lemma: critical shperes don't accumulate}, \ref{lemma: getting epsilon small enough} and \ref{lemma: retraction of B_r onto B_r^s}.

\begin{corollary} \label{corol: descrete sequence of expanding simplicial balls}
    For any $r>0$ such that  $\mathbf S_{r,v}$ is a critical level sphere, there is an $0< \epsilon < r-\sqrt{r^2 - \frac{1}{4}} $, such that every vertex $x\subset \cal C$, where $d(x,v)<r$, is contained in  $\mathbf B_{r-\epsilon, v}^s$.
\end{corollary}

Along with the geodesic extension property,  this corollary tells us that there is an infinite sequence of engulfing simplicial balls in $\cal C$.


\section{Expanding sequence in $\tld M$}

All that remains to prove theorem \ref{theorem: universal cover is $R^3$} is to show that the interior of the sub-manifold in $\tld M$ that corresponds to a simplicial ball in $\cal C$ is homeomorphic to $\R[3]$.  We do this by defining a process to get from $\mathbf B_{r-\epsilon,v}^s$ to $\mathbf B_{r,v}$, such that at each step the interior of the corresponding sub-manifold in $\tld M$ is an open ball.


\begin{definition}
    For a simplicial set $C$ and a vertex $x\subset C$, let  $adj_C(x)=st(x)\cap C$, that is the set of simplices in $C$ adjacent to $x$, and $adj_C^o(x)=st^o(x)\cap C$.
\end{definition}

Note we can get $adj_C^o(x)$ from $ adj_C(x)$ by removing its intersection with $lk(x)$.  Also by lemma \ref{lemma: getting epsilon small enough} we know that for some vertex $x\in V_{r,v}$ and $0<\epsilon< r-\sqrt{r^2 - \frac{1}{4}}$, the level sphere $\mathbf S_{r-\epsilon,v}$ intersects every edge of $adj_{\mathbf B_{r,v}}^o(x)$ transversely.


\begin{lemma}
    For $x\in V_{r,v}$, $\adj(x)$ is convex.
\end{lemma}

\begin{remark}
The idea behind this proof is the same as that used for the proof of lemma \ref{lemma: st(v) is convex}.
\end{remark}

\begin{proof}
Let  $\Gamma = lk(x)\cap \adj(x)$ and $\Pi$ be the set of free edges of $\adj(x)$ not in $\Gamma$,  that is the free edges  adjacent to $x$. Let $D=\mathbf B_{r,v}\cap st(x)$.  By convexity of $\mathbf B_{r,v}$ and $st(x)$, $D$ is convex.  Also note that $\adj(x)\subset D$.  Let $C=D-\adj(x)$ and $C_i$'s be the components of $C$.  As all the vertices in $D$ are in $\adj(x)$, then $\adj(x)\cap \cal C^0 = \emptyset$.  Also $C\cap \cal C^1\subset lk(x)$ and each edge in the $lk(x)$ is free in $st(x)$.  Therefore each $C_i$ is contained in a single face of $st(x)$.  For any $C_i$, $\overline{C_i}\cap \adj(x)$ is some edge $E_i\subset \Pi$.  Therefore if $\adj(x)$ is not convex then there must be two points $\{x,y\}\subset \adj(x)$ such that the image, $\alpha$, of the geodesic joining them is  not entirely contained in $\adj(x)$.  However $D$ is convex so $\alpha\subset D$.  This means that for at least one $C_i$, that $C_i\cap \alpha \not= \emptyset$.  Moreover, this means that $\alpha $ must intersect $E_i$ at least twice, which contradicts   the unique geodesic property of CAT(0) spaces. Thus $\adj(x)$ is convex.
\end{proof}


\begin{lemma}\label{lemma: max 2 vert of  face in V_rv}
    Let $F$ be a face of $\cal C$ then all three vertices of $F$ cannot be in $V_{r,v}$.
\end{lemma}

\begin{proof}
By lemma \ref{lemma: getting epsilon small enough} we know that if F's three vertices are in $V_{r,v}$ then there is an $\epsilon > 0$ such that $F\cap \mathbf S_{r-\epsilon,v}$ has three components. However by lemma \ref{corollary: S_rv intersect F can't have three components} we know this is not possible.
\end{proof}

\begin{lemma}\label{lemma: E_xx' is a free edge in adj(x)}
    If $x$ and $x'$ are  in $V_{r,v}$ and neighbouring vertices, then the edge between them is free in $\adj(x)$.
\end{lemma}

\begin{proof}
Let $E$ be the edge between $x$ and $x'$.  Assume that the lemma is not true and that $E$ is adjacent to more than one face in $\adj(x)$.  By lemma \ref{lemma: max 2 vert of  face in V_rv} we know that the third vertex in any such face is contained in $\mathbf B_{r-\epsilon,v}$.  However as the metric ball is retractible in $\cal C$ we know that for all but one of these  faces in $\adj(x)$ adjacent to $E$ that their closest point to $X$ is in $int(E)$.  Let $F$ be such a face, that is $y=\{z\in F: d_F(v) = d(z,v)\}\subset int(E)$.  By lemma \ref{lemma: B&H convex complete sets} we know that $y$ is a single point.  Let $x''$ be the third vertex of $F$. Then $d(x'',v)<r$.  Let $E_1$ be the edge between $x$ and $x''$ and let $E_2$ be the edge between $x'$ and $x''$.  Also let $y_i=\{z\in E_i:d_{E_i}(v)=d(z,v)\}$.  As $F\subset \mathbf B_{r,v}$ and $d(x,v)=d(x',v)=r$ we know that $y_1\not= x$ and $y_2\not= x'$.  The closest point of $F$ to $v$ is $int(E)$ and for $d_F(v)<t<r$ $S_{t,v}\cap F$ is a set of circular arcs. Then it is not possible for $y_1=y_2= x''$ and at least one $y_i$ must be in $int(E_i)$.  Let it be $E_1$.   If $y_2=x''$ then $S_{t,v}\cap F$, for $d_{E_1}(v)<t<d(x'',v)$, has two components connecting the same pair of edges. However by lemma \ref{lemma: S_r can't have parallel components in F} we know that this is not possible.  Therefore $y_2\subset int(E_2)$.  Let $d_{E_1}(v)\geq d_{E_2}(v)$, then  $S_{t,v}\cap F$, for $d_{E_1}(v)<t<d(x'',v)$ has three components. However by lemma \ref{corollary: S_rv intersect F can't have three components} we know this is not possible.
\end{proof}

To produce the process for expanding from one simplicial ball to the next we need to define a  simplicial join.

\begin{figure}[h]
  \begin{center}
      \includegraphics[width=6cm]{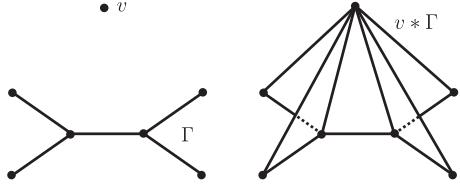}
   \end{center}
  \caption{Simplicial cone.}
  \label{fig: Simplicial join.}
\end{figure}

\begin{definition}
    Let $v$ be a vertex of $\cal C$ and $\Gamma$ be a graph in $lk(v)$, then their \textbf{simplicial cone},  $a\ast \Gamma$, is the simplicial subset of $st(v)$ produced by coning $\Gamma$ to $v$, as shown in figure \ref{fig: Simplicial join.}.
\end{definition}

Note that in combinatorial topology this process is more generally known as a simplicial join. If $\Gamma = lk(x)$, then $st(x) = v \ast \Gamma$.

\begin{definition}
    Let $x\in \cal C^0$ and $T\subset lk(x)$ is a connected tree, then using the metric on $lk(x)$, as defined in the proof of lemma \ref{lemma: C has non+'ve curvature}, the \textbf{diameter} of  $T$ is $diam(T)= max\{d(a,b):a,b\in T\}$.
\end{definition}

For some simplicial set $C \subseteq \cal C$ let $\tld C$ be the corresponding sub-manifold of $\tld M$. The following lemma gives the `move' used to expand a simplicial set in $\cal C$ so that the interior of the corresponding set in $\tld M$ remains homeomorphic to $\R[3]$.


\begin{lemma}\label{lemma: manifold corresponding to join}
    Let $x\in \cal C^0$ and $T\subset lk(x)$ be a connected tree such that $diam(T)\leq \pi$, then the interior of the manifold in $\tld M$ corresponding to $x\ast T$ is homeomorphic to $\R[3]$.
\end{lemma}

\begin{proof}
Let $F= int(\tld T \cap \tld x)\subset \tld M$.  By lemmas \ref{lemma: interior of the union of two balls is an open ball} and \ref{lemma: condition for product structure} it is sufficient to show that $int(\tld T) \homeo \R[3]$ and $F \homeo \R[2]$.  Let $E$ be an edge of $T$ and $a,b\subset E$ the two vertices of $E$.  Then $\tld E$ is two lifts of handlebodies glued along a face $G=\tld a \cap \tld b \subset \tld X$.  As $a$ and $b$ are joined by single edge in $E$, $int(G) \homeo \R[2]$.  By  lemmas \ref{lemma: interior of the union of two balls is an open ball} and \ref{lemma: condition for product structure}, $\tld E\homeo \R[3]$.  Let $T', T''\subseteq T $ be  connected trees such that $T''$ is produced from $T'$ by adding a single edge.  Let $a$ be the vertex contained in $T''$ but not in $T'$.   As $a$ is joined to $T'$ by a single edge, $int(\tld a\cap \tld T')\homeo \R[2]$. Therefore if $\tld T' \homeo \R[3]$, then $\tld T'' \homeo \R[3]$. Thus by using induction $\tld T \homeo \R[3]$.  Taking the embedded view of $\cal C$,  $\bndry \tld x \cap (\tld{x\ast T})$ is a graph dual to $F$ and homeomorphic to $T$.  As $diam(T)\leq \pi$ by the disk-condition $F$ is a missing boundary disk in $\bndry \tld x$ and thus $int(F) \homeo \R[2]$ .
\end{proof}

The following corollary is a generalisation of the previous lemma.

\begin{corollary}\label{lemma: simplicial join gives expanding Ball}
Let $C\subset \cal C$ be a connected simplicial set such that $int(\tld C)\homeo \R[3]$  and $x\in \cal C^0$ such that $T = lk(x)\cap C$ is a connected tree and $diam(T)< \pi$, then the interior of the manifold corresponding to $C\cup(x\ast T)$ is homeomorphic to $\R[3]$.
\end{corollary}

This comes from the observation that $\tld x \cap \tld C = \tld x \cap \tld T$, where $\tld x$, $\tld C$ and $\tld T$ are the manifolds in $\tld M$ corresponding to $x$, $C$ and $T$ respectively.


Now the next step is to show that the intersection between the link of a vertex in $V_{r,v}$ and the simplicial ball $B_{r,v}^s$ is in fact a tree. Note that $adj_{\mathbf B_{r,v}^s} (x) = x \ast T$.

\begin{lemma} \label{lemma:adj(x)=x*T}
If $x\in V_{r,v}$, then $adj_{\mathbf B_{r,v}^s} (x) \cap lk(x) = T$ is a connected tree such that $diam(T)\leq \pi$.
\end{lemma}

\begin{proof}
By lemma \ref{lemma: critical shperes don't accumulate} we can find an $0<\epsilon < r-\sqrt{r^2 - \frac{1}{4}}$ such that the level sphere $\mathbf S_{r- \epsilon,v}$ intersects every edge in $\adj(x)$ adjacent to $x$, at least once transversely and $B_{r-\epsilon,v}^s$ contains all vertices that are distance less than $r$ from $v$. Let $\Gamma = \mathbf S_{r-\epsilon} \cap \adj(x)$.  Then $\Gamma = \{\Gamma_i\}$ is a set of graphs in $\adj(x)$.

First we want to show $\Gamma$ has a single component that separates $x$ from all the other vertices of $\adj(x)$. As $\mathbf S_{r-\epsilon,v}$ is separating,  $\Gamma$ is separating in $\adj(x)$.  However, as $\adj(x)$ is convex and thus simply connected, each component of $\Gamma$ must be separating.   Therefore as $\mathbf S_{r-\epsilon,v}$ intersects every edge of $\adj(x)$ adjacent to $x$, that $\Gamma$ separates $x$ from all other vertices in $\adj(x)$.

Let $\Gamma'$ be a subset of components of $\Gamma$ with a minimal number of components that separate $x$ from all other vertices in $\adj(x)$.  We want to show that $\Gamma'$ is connected, that is contains only one component.  If $\Gamma'$ is not connected then, as $x$ is adjacent to every face in $\adj(x)$, there must be a face $F$ of $\adj(x)$, such that $\Gamma' \cap F$ is two embedded arcs.  By lemma  \ref{corollary: S_rv intersect F can't have three components} $\Gamma' \cap F$ cannot contain more than two. For two neighbouring vertices, $a$ and $b$, let $E_{ab}$ be the edge joining them.   By Lemma \ref{lemma: max 2 vert of  face in V_rv} we know that for any face of $\adj(x)$ that at most two of its vertices are in $V_{r,v}$.  First let $F$ be a face such that two of its vertices, $a_1$ and $a_2$, are in $\mathbf B_{r-\epsilon,v}$.   Then the edge $E_{a_1a_2}$ must be contained in $\mathbf B_{r-\epsilon,v}$ and thus $\Gamma$ does not intersect it. As $a_i$'s are contained in $\mathbf B_{r-\epsilon,v}$ and $x$ is not, $\Gamma$ intersects $E_{xa_i}$ exactly once.  Therefore $F\cap \Gamma$ contains only one arc.

The second case is when only one vertex, $a$, of $F$ is contained in $\mathbf B_{r-\epsilon,v}$.  Therefore the third vertex $x'$ must be in $V_{r,v}$.  Therefore the edges $E_{xa}$ and $E_{x'a}$ both intersect $\Gamma$ once. As both $x$ and $x'$ are not contained in $\mathbf B_{r-\epsilon,v}$, by definition of $\epsilon$ we know that  $E_{xx'}$ intersects $\mathbf B_{r-\epsilon,v}$ in its interior.  Thus $E_{xx'}$ intersects $\Gamma$ twice.  This means that $\Gamma \cap F$ is two properly embedded arcs in $F$, one between $E_{xx'}$ and $E_{xa}$ and the other between $E_{xx'}$ and $E_{x'a}$.

By lemma \ref{lemma: E_xx' is a free edge in adj(x)} we know that $E_{xx'}$ is a free edge of $\adj(x)$. Thus the arc of $\Gamma\cap F$ that runs between $E_{xx'}$ and $E_{x'a}$ is a component of $\Gamma$ as it runs between free edges.  This means that $\Gamma'$ will either contain the edge between $E_{xx'}$ and $E_{x'a}$ or the component that contains the other arc of $\Gamma \cap F$.  Therefore every face of $\adj(x)$ contains only one edge of $\Gamma'$.  Thus $\Gamma'$ is connected.

Next we want to show that $\Gamma'$ is a tree.  As both $\adj(x)$ and $\mathbf B_{r-\epsilon,v}$ are convex,  $\adj(x) \cap \mathbf B_{r-\epsilon,v}$ is convex and does not contain $x$.  Therefore if $\Gamma'$ contained a loop then $\adj(x) \cap \mathbf B_{r-\epsilon,v}$ must contain the disk bounded by that loop.  However by the cone structure of $ st(x)$ such a disk must contain $x$.   Thus $\Gamma'$ is a tree.

When we cut $\adj(x)$ along $\Gamma'$, because of the cone structure of $st(x)$,   the component that doesn't contain $x$ has the product structure $\Gamma' \cross I$.  Therefore $T\homeo  \Gamma'$.  As $T\subset B_{r,v}\cap lk(x)$ and by convexity of $B_{r,v}$, $diam(T) \leq \pi$.
\end{proof}

\begin{lemma}\label{lemma: lk(x) intersect B_(r-e) is a connected tree}
If $x\in V_{r,v}$, then $T' = \adj(x)\cap \mathbf B_{r-\epsilon, v}^s\subset lk(x)$ is a connected tree and $diam(T')\leq\pi$.
\end{lemma}

\begin{proof}
Let $T= \adj(x) \cap lk(x)$.  Clearly $T' \subset T $, therefore $diam(T') \leq \pi$ and by lemma \ref{lemma:adj(x)=x*T} each component of $T'$ is  a tree, so all we need to prove is that $T'$ is connected. We can get $T'$ from $T$ by removing all the vertices $V_{r,v}\cap T$ and their adjacent edges.  However we know from lemma  \ref{lemma: E_xx' is a free edge in adj(x)} that the edge between any vertex $V_{r,v}\cap T$ and x is free in $\adj(x)$,   therefore any vertex in $V_{r,v}\cap T$ has degree one in $T$. This means that all the vertices in $V_{r,v}\cap T$ are not separating in $T$ and thus $T'$ is connected.
\end{proof}


\begin{lemma}
If the interior of the manifold in $\tld M$ corresponding to  $\mathbf B_{r-\epsilon,v}^s$ is an open ball, then the interior of the manifold corresponding to  $\mathbf B_{r,v}^s$ is an open ball.
\end{lemma}

\begin{proof}
To prove the above lemma we will define a process for getting from $\mathbf B_{r-\epsilon,v}^s$ to $\mathbf B_{r,v}^s$, then show that at each step the simplicial set corresponds to a missing boundary ball in $\tld M$ and then that we do indeed get $\mathbf B_{r,v}^s$.

Let $V_{r,v} = \{x_1, x_2, ...\}$ and $C_0 = \mathbf B_{r-\epsilon,v}$.  Then let $\Gamma_i = lk(x)\cap C_{i-1}$ and $C_i = C_{i-1}\cup (\Gamma_i \ast x_i)$.  To show that each $C_i$ corresponds to a missing boundary ball, by lemma \ref{lemma: simplicial join gives expanding Ball} we just need to show that each $\Gamma_i$ is a connected tree.

We know that $(\mathbf B_{r-\epsilon,v}^s \cap lk(x_i))\subseteq \Gamma_i \subseteq (\mathbf B_{r,v}^s \cap lk(x_i))$ and by lemma \ref{lemma:adj(x)=x*T} and lemma \ref{lemma: lk(x) intersect B_(r-e) is a connected tree} both $\mathbf B_{r-\epsilon,v}^s \cap lk(x_i)$ and $\mathbf B_{r,v}^s \cap lk(x_i)$ are connected trees. Also by lemma \ref{lemma: E_xx' is a free edge in adj(x)} any vertices of $V_{r,v}$ in $\mathbf B_{r,v}^s \cap lk(x_i)$ are not separating.  Therefore each $\Gamma_i$ is a connected tree.

As $\Gamma_i\subseteq (\mathbf B_{r,v}^s \cap lk(x_i))$, then $x_i \ast \Gamma_i \subseteq \adj(x)$. Therefore we know that $C_i\subseteq B_{r,v}^s$.  If $C_i \subset \mathbf B_{r,v}^s$ then $\mathbf B_{r,v}^s - C_i$ must contain at least one vertex $x$, and $x$ must be in $V_{r,v}$.  Therefore once all the vertices of $V_{r,v}$ have been added to $\mathbf B_{r-\epsilon,v}$ by this process, the resulting simplicial set is $\mathbf B_{r,v}^s$
\end{proof}

\begin{proof} \emph{(of theorem \ref{theorem: universal cover is $R^3$})}
To prove lemma \ref{theorem: universal cover is $R^3$} we use the previous lemma to  produce an infinite sequence of expanding open balls and by Brown \cite{br}, this implies that $\tld M$ is homeomorphic to $\R$.  By lemma \ref{corol: descrete sequence of expanding simplicial balls}, for any vertex $v\in \cal C$ there is a sequence $R=\{r_i\}$ such that $\mathbf B_{r_i,v}^s \subset  \mathbf B_{r_{i+1},v}^s$ and by the previous lemma we know that the interior of the manifolds corresponding to the sequence $\{\mathbf B_{r_i,v}^s\}$ of simplicial balls are open balls in $\tld M$. Thus we have an infinite sequence of engulfing open balls and the universal cover is homeomorphic to $\R[3]$.
\end{proof}

\bibliographystyle{plain}
\bibliography{bibfile}

\end{document}